\documentclass[10pt]{article}
\usepackage{lmodern}
\usepackage{amsmath}
\usepackage[T1]{fontenc}
\usepackage[utf8]{inputenc}
\usepackage{authblk}
\usepackage{amsfonts}
\usepackage{graphicx}
\usepackage{rotating}
\usepackage{amssymb}
\usepackage[english]{babel}
\usepackage{color}
\usepackage{amsthm}
\usepackage{mathrsfs}
\usepackage{makecell}
\usepackage{microtype}
\usepackage{enumerate}
\usepackage{mathscinet}
\usepackage{array}
\usepackage{multirow}
\usepackage{booktabs}
\usepackage{enumerate}
\usepackage[cal=boondoxo,bb=ams]{mathalfa}
\usepackage{hyperref}
\hypersetup{hidelinks}
\usepackage{titlesec}
\newtheorem{theorem}{Theorem}[section]
\newtheorem{prop}{Proposition}[section]

\newtheorem{remark}{Remark}[section]

\newcommand{\ml}{\mathcal}
\newcommand{\mb}{\mathbb}

\DeclareMathOperator{\intt}{int}
\DeclareMathOperator{\extt}{ext}
\DeclareMathOperator{\bdd}{bdd}
\DeclareMathOperator{\divv}{div}

\title{Some remarks on large-time behaviors for the linearized compressible Navier-Stokes equations}
\author[1]{Wenhui Chen\thanks{Wenhui Chen (wenhui.chen.math@gmail.com)}}
\affil[1]{School of Mathematics and Information Science, Guangzhou University, 510006 Guangzhou, China}
\author[2]{Ryo Ikehata\thanks{Ryo Ikehata (ikehatar@hiroshima-u.ac.jp)}}
\affil[2]{Department of Mathematics, Division of Educational Sciences, Graduate School of Humanities and Social Sciences, Hiroshima University, 739-8524 Higashi-Hiroshima, Japan}

\date{}
\begin{document}

\maketitle
\begin{abstract}
	\medskip
In this paper, we consider the linearized compressible Navier-Stokes equations in the whole space $\mathbb{R}^n$. Concerning initial  datum with suitable regularities, we introduce a new threshold $|\mb{B}_0|=0$ to distinguish different large-time behaviors. Particularly in the lower-dimensions, optimal growth estimates ($n=1$ polynomial growth, $n=2$ logarithmic growth) hold when $|\mb{B}_0|>0$, whereas optimal decay estimates hold when $|\mb{B}_0|=0$. Furthermore, we derive asymptotic profiles of solutions with weighted $L^1$ datum as large-time.\\
	
	
	\noindent\textbf{AMS Classification (2020)}  35Q35, 35B40, 76N10, 35Q30, 
	
\end{abstract}
\fontsize{12}{15}
\selectfont
\section{Introduction}
We are interested in exploring large-time asymptotic behaviors for the linearized compressible Navier-Stokes equations in the whole space $\mb{R}^n$ with any $n\geqslant 1$, namely,
\begin{align}\label{Eq_Linearized_NS}
\begin{cases}
\rho_t+\gamma\divv v=0,&x\in\mb{R}^n,\ t>0,\\
v_t-\alpha\Delta v-\beta\nabla\divv v+\gamma\nabla\rho=0,&x\in\mb{R}^n,\ t>0,\\
\rho(0,x)=\rho_0(x),\ v(0,x)=v_0(x),&x\in\mb{R}^n,
\end{cases}
\end{align}
with the density $\rho=\rho(t,x)\in\mb{R}$ and the velocity $v=v(t,x)\in\mb{R}^n$, where $\alpha,\gamma$ are positive constants and $\beta$ is a non-negative constant. For briefness, we denote $v^{(k)}=v^{(k)}(t,x)\in\mb{R}$ to be the $k$-th element of the vector $v(t,x)$ with $k=1,\dots,n$. Our main purpose of this work is to investigate a new threshold with the following quantity:
\begin{align}\label{Intro_B0}
	|\mb{B}_0|^2:=|P_{|D|\rho_0}|^2+|P_{|D|v_0}|^2
\end{align}
to distinguish two different large-time asymptotic behaviors of solutions, where we denote
\begin{align*}
	P_{|D|\rho_0}=\int_{\mb{R}^n}|D|\rho_0(x)\mathrm{d}x\ \ \mbox{and}\ \ P_{|D|v_0}=\left(\int_{\mb{R}^n}|D|v_0^{(1)}(x)\mathrm{d}x,\cdots,\int_{\mb{R}^n}|D|v_0^{(n)}(x)\mathrm{d}x\right)
\end{align*}
with the pseudo-differential operator $|D|$ owning its symbol $|\xi|$. Especially, the solutions $(\rho,v)$ grow polynomially when $n=1$ and logarithmically when $n=2$ in the $L^2$ framework if $|\mb{B}_0|>0$, whereas the solutions decay polynomially if $|\mb{B}_0|=0$. For the higher-dimensions $n\geqslant 3$, the threshold for \eqref{Intro_B0} distinguishes two kinds of optimal decay estimates with different decay rates.

Let us now recall several historical background for the compressible Navier-Stokes equations in $\mb{R}^n$. It describes the motion of the general isotropic Newtonian fluids. In 1979, the authors of \cite{Mat-Nishi=1979} demonstrated global (in time) existence of solutions to the three-dimensional case with small datum in $H^3\cap L^1$, where the velocity fulfills heat-type decay estimates. Here, the nomenclature \emph{heat-type decay estimates} means $L^2$ estimates of solutions with decay rates $(1+t)^{-\frac{n}{4}}$ for $n\geqslant 1$. The solution (in time) asymptotic to the one of its linearized problem was shown in \cite{Kawa-Mat-Nishi=1979}. Later, the authors of \cite{Hoff-Zumbrun=1995} derived some detailed large-time behaviors of solutions in the $L^p$ framework with $1\leqslant p\leqslant \infty$, and the second-order asymptotic profile was obtained in \cite{Kagei-Okita=2017}. We mention that the diffusion-wave structure of the linearized compressible Navier-Stokes equations \eqref{Eq_Linearized_NS} have been deeply studied in \cite{Kobayashi-Shibata=2002,Charao-Ikehata=2019}, in which the solutions $(\rho,v)$ to the Cauchy problem \eqref{Eq_Linearized_NS} fulfill general heat-type decay estimates.  Up to our best knowledge, the heat-type decay properties are crucial in studies of the linearized or nonlinear compressible Navier-Stokes equations. 

In the present paper, we consider the Cauchy problem \eqref{Eq_Linearized_NS} with initial datum $(|D|\rho_0,|D|v_0)$ belonging to some (weighted) $L^1$ spaces. By employing the WKB analysis and the Fourier analysis, we derive optimal estimates and asymptotic profiles of the solutions $(\rho,v)$ under different assumptions on initial datum, which indicates a new threshold $|\mb{B}_0|=0$ (defined in \eqref{Intro_B0} before) for heat-type decay properties and growth properties as large-time. Totally, the main results with explanations  will be stated in Section \ref{Section-Main-result}, and their proofs will be given in Section \ref{Section-Fourier-space}. Because our philosophy can be applied in some linearized models in compressible fluids, we will propose some remarks in Section \ref{Section_Final_Remark} to end this paper.

\medskip
\noindent\textbf{Notations:} The relation $f\lesssim g$ means that there exists a positive constant $C$ fulfilling $f\leqslant Cg$, which may be changed in different lines, analogously, for $f\gtrsim g$. Furthermore, the asymptotic relation  $f\simeq  g$ holds if and only if $f\lesssim g$ and $f\gtrsim g$ simultaneously, which will be used frequently in optimal estimates. The notation $\circ$ denotes the inner product in the Euclidean space. We denote the homogeneous Bessel space $\dot{H}^{1}_1:=|D|^{-1}L^1$. Let us recall the weighted $L^1$ spaces with $s>0$ as follows:
\begin{align*}
	L^{1,s}:=\left\{f\in L^1  :\ \|f\|_{L^{1,s}}:=\int_{\mb{R}^n}(1+|x|)^s|f(x)|\mathrm{d}x<\infty \right\}.
\end{align*}
The  means of a summable function $f$ are denoted by $P_f:=\int_{\mb{R}^n}f(x)\mathrm{d}x$ and $Q_f:=\int_{\mb{R}^n}xf(x)\mathrm{d}x$, whose form is $Q_f=(Q_f^{(1)},\dots,Q_f^{(n)})$. We introduce the following time-dependent function:
\begin{align}\label{Decay-fun}
	\ml{D}_n(t):=\begin{cases}
		\sqrt{t}&\mbox{when}\ \ n=1,\\
		\sqrt{\ln t}&\mbox{when}\ \ n=2,\\
		t^{\frac{1}{2}-\frac{n}{4}}&\mbox{when}\ \ n\geqslant 3,
	\end{cases}
\end{align}
to complete this introduction.

\section{Main results}\label{Section-Main-result}
We state the first result under the assumption $|\mb{B}_0|>0$. It implies some growth properties of solutions in the lower-dimensions.
\begin{theorem}\label{Thm-Optimal-B0}
Let us consider the linearized compressible Navier-Stokes equations \eqref{Eq_Linearized_NS} with initial datum $(\rho_0,v_0)\in (L^2\cap \dot{H}^1_1)^{1+n}$ and $|\mb{B}_0|>0$. Then, the solutions fulfill the optimal estimates
\begin{itemize}
	\item for the lower-dimensions $n=1,2$:
	\begin{align*}
	\ml{D}_n(t)|\mb{B}_0|\lesssim\left\|(\rho,v)(t,\cdot)\right\|_{(L^2)^{1+n}}\lesssim\ml{D}_n(t)\|(\rho_0,v_0)\|_{(L^2\cap\dot{H}^1_1)^{1+n}};
	\end{align*}
	\item for the higher-dimensions $n\geqslant 3$:
	\begin{align*}
	\ml{D}_n(t)|\mb{B}_0|\lesssim\left\|\left(\rho-\mathrm{e}^{\frac{\alpha+\beta}{2}\Delta t}\rho_0,v-\mathrm{e}^{\alpha\Delta t}v_0\right)(t,\cdot)\right\|_{(L^2)^{1+n}}\lesssim\ml{D}_n(t)\|(\rho_0,v_0)\|_{(L^2\cap\dot{H}^1_1)^{1+n}};
	\end{align*}
\end{itemize}
as large-time $t\gg1$, where the time-dependent function $\ml{D}_n(t)$ was defined in \eqref{Decay-fun}. Furthermore, by assuming $(|D|\rho_0,|D|v_0)\in (L^{1,1})^{1+n}$ additionally, the solutions fulfill the refined estimates
\begin{align}
&\left\|\left(\rho-\mathrm{e}^{\frac{\alpha+\beta}{2}\Delta t}\rho_0-\widetilde{\rho},v-\mathrm{e}^{\alpha\Delta t}v_0-\widetilde{v}\right)(t,\cdot)\right\|_{L^2}\lesssim t^{-\frac{n}{4}}\left(\|(|D|\rho_0,|D|v_0)\|_{(L^{1,1})^{1+n}}+\|(\rho_0,v_0)\|_{(L^2)^{1+n}}\right),\label{Est-002}
\end{align}
for large-time $t\gg1$, where the corresponding profiles are described by
\begin{align*}
\widetilde{\rho}(t,x)&:=J_0(t,x)P_{|D|\rho_0}+J_1(t,x)\circ P_{|D|v_0},\\
\widetilde{v}(t,x)&:=J_1(t,x)P_{|D|\rho_0}+J_2(t,x)\ml{F}^{-1}\left(\frac{\xi(\xi\circ P_{|D|v_0})}{|\xi|^2}\right),
\end{align*}
with the auxiliary functions defined via the Fourier multipliers
\begin{align*}
J_0(t,x)&:=\ml{F}^{-1}_{\xi\to x}\left(\frac{1}{|\xi|}\big(\cos(\gamma|\xi|t)-1\big)\mathrm{e}^{-\frac{\alpha+\beta}{2}|\xi|^2t}\right),\\ J_1(t,x)&:=-i\,\ml{F}^{-1}_{\xi\to x}\left(\frac{\sin(\gamma|\xi|t)}{|\xi|}\mathrm{e}^{-\frac{\alpha+\beta}{2}|\xi|^2t}\frac{\xi}{|\xi|}\right),\\
J_2(t,x)&:=\ml{F}^{-1}_{\xi\to x}\left(\frac{1}{|\xi|}\left(\cos(\gamma|\xi|t)\mathrm{e}^{-\frac{\alpha+\beta}{2}|\xi|^2t}-\mathrm{e}^{-\alpha|\xi|^2t}\right)\right).
\end{align*}
\end{theorem}
\begin{remark}By assuming $(\rho_0,v_0)\in(L^2\cap\dot{H}^1_1)^{1+n}$ associated with $|\mb{B}_0|>0$, the unknowns $(\rho,v)$ in the $L^2$ framework, grow polynomially when $n=1$; logarithmically when $n=2$. Those estimates are optimal for large-time. This phenomenon seems new for the mathematical models in compressible fluids.
\end{remark}
\begin{remark}
	With additionally weighted $L^1$ datum when $|\mb{B}_0|>0$, by subtracting the profiles $\widetilde{\rho}(t,\cdot)$ and $\widetilde{v}(t,\cdot)$ in the $L^2$ norm, the obtained estimates have been improved by $t^{-\frac{3}{4}}$ when $n=1$, $(t\log t)^{-\frac{1}{2}}$ when $n=2$, and $t^{-\frac{1}{2}}$ when $n\geqslant 3$. These new large-time profiles have the diffusion-wave structures
	\begin{align*}
	\ml{F}^{-1}_{\xi\to x}\left(\frac{1}{|\xi|}\sin(c_0|\xi|t)\mathrm{e}^{-c_1|\xi|^2t}\right)\ \ \mbox{as well as}\ \ \ml{F}^{-1}_{\xi\to x}\left(\frac{1}{|\xi|}\left(\mathrm{e}^{-c_2|\xi|^2t}-\cos(c_0|\xi|t)\mathrm{e}^{-c_1|\xi|^2t}\right)\right)
	\end{align*}
	with $c_0,c_1,c_2>0$, which consist of the singularity $|\xi|^{-1}$ for small frequencies, the oscillations $\sin(c_0|\xi|t)$ and $\cos(c_0|\xi|t)$, the diffusion $\mathrm{e}^{-c_1|\xi|^2t}$. 
\end{remark}
\begin{remark}
	Let us compare the classical results with $L^2\cap L^1$ datum for the linearized compressible Navier-Stokes equations (see \cite{Hoff-Zumbrun=1995,Kobayashi-Shibata=2002,Ikehata-Onodera=2017,Charao-Ikehata=2019} and references therein). By assuming $(\rho_0,v_0)\in (L^2\cap L^{1})^{1+n}$ as well as $|P_{\rho_0}|+|P_{v_0}|\neq0$, the optimal heat-type decay estimates $\|(\rho,v)(t,\cdot)\|_{(L^2)^{1+n}}\simeq t^{-\frac{n}{4}}$ hold 
	for any $n\geqslant 1$ and $t\gg1$. Nevertheless, by changing the datum spaces into $(L^2\cap \dot{H}^1_1)^{1+n}$ and taking $|\mb{B}_0|>0$, some new diffusion-wave behaviors and optimal growth properties $\|(\rho,v)(t,\cdot)\|_{(L^2)^{1+n}}\simeq \ml{D}_n(t)$ when $n=1,2,$ occur.
\end{remark}

Let us turn to the second result with $|\mb{B}_0|=0$ but $|\mb{B}_1|\neq0$, where we introduce another quantity
\begin{align*}
|\mb{B}_1|^2:=|Q_{|D|\rho_0}|^2+\sum\limits_{k=1}^n\left(|Q_{|D|v_0^{(k)}}|^2+|Q_{|D|v_0^{(k)}}^{(k)}|^2 \right).
\end{align*}
Then, it recovers heat-type decay properties of solutions even with initial datum $(|D|\rho_0,|D|v_0)$ in some weighted spaces. In other words, the condition $|\mb{B}_0|=0$ is the crucial threshold to determine two different large-time behaviors of the solutions.
\begin{theorem}\label{Thm_Optimal-B1}
	Let us consider the linearized compressible Navier-Stokes equations \eqref{Eq_Linearized_NS} with initial datum $(|D|\rho_0,|D|v_0)\in (L^{1,2})^{1+n}$, $(\rho_0,v_0)\in (L^2)^{1+n}$ and $|\mb{B}_0|=0$. Then, the solutions fulfill the optimal estimates
		\begin{align*}
		t^{-\frac{n}{4}}|\mb{B}_1|\lesssim\left\|(\rho,v)(t,\cdot)\right\|_{(L^2)^{1+n}}\lesssim t^{-\frac{n}{4}}\left(\|(|D|\rho_0,|D|v_0)\|_{(L^{1,1})^{1+n}}+\|(\rho_0,v_0)\|_{(L^2)^{1+n}}\right)
		\end{align*}
	for any $n\geqslant 1$, as large-time $t\gg1$. Furthermore,   the solutions fulfill the refined estimates
	\begin{align*}
	&\left\|\left(\rho-\widetilde{\rho}^0,v-\widetilde{v}^0\right)(t,\cdot)\right\|_{(L^2)^{1+n}}\lesssim t^{-\frac{1}{2}-\frac{n}{4}}\left(\|(|D|\rho_0,|D|v_0)\|_{(L^{1,2})^{1+n}}+\|(\rho_0,v_0)\|_{(L^2)^{1+n}}\right),
	\end{align*}
	for large-time $t\gg1$, where the corresponding profiles are described by
	\begin{align*}
	\widetilde{\rho}^0(t,x)&:=i\ml{F}^{-1}_{\xi\to x}\left(\frac{\cos(\gamma|\xi|t)}{|\xi|}\mathrm{e}^{-\frac{\alpha+\beta}{2}|\xi|^2t}\xi\right)\circ Q_{|D|\rho_0}+i\sum\limits_{k=1}^n\ml{F}^{-1}_{\xi\to x}\left(\widehat{J}_1^{(k)}\xi\right)\circ Q_{|D|v_0^{(k)}},\\
	\widetilde{v}^0(t,x)&:=i\ml{F}^{-1}_{\xi\to x}\left(\widehat{J}_1\xi\right)\circ Q_{|D|\rho_0}+i\ml{F}^{-1}_{\xi\to x}\left(\frac{\xi}{|\xi|}\mathrm{e}^{-\alpha|\xi|^2t}\right)\circ Q_{|D|v_0}+i\sum\limits_{k=1}^n\ml{F}^{-1}_{\xi\to x}\left(\frac{\xi}{|\xi|^2}\widehat{J}_2\xi_k\xi\right)\circ Q_{|D|v_0^{(k)}}.
	\end{align*}
\end{theorem}
\begin{remark}
When $|D|\rho_0$ and $|D|v_0^{(k)}$ are odd functions with respect to $x_j$ for all $k=1,\dots,n$, we immediately obtain $|P_{|D|\rho_0}|=|P_{|D|v_0^{(k)}}|=0$, in other words, $|\mb{B}_0|=0$.
\end{remark}
\begin{remark}
Summarizing Theorems \ref{Thm-Optimal-B0} and \ref{Thm_Optimal-B1}, we may understand the importance of $|\mb{B}_0|$ by the next way:
\begin{align*}
\|(\rho,v)(t,\cdot)\|_{(L^2)^{1+1}}&\simeq\begin{cases}
\sqrt{t}&\mbox{when}\ \ |\mb{B}_0|>0,\\
t^{-\frac{1}{4}}&\mbox{when}\ \ |\mb{B}_0|=0,
\end{cases}\\
\|(\rho,v)(t,\cdot)\|_{(L^2)^{1+2}}&\simeq\begin{cases}
\sqrt{\ln t}&\mbox{when}\ \ |\mb{B}_0|>0,\\
t^{-\frac{1}{2}}&\mbox{when}\ \ |\mb{B}_0|=0,
\end{cases}
\end{align*}
for large-time $t\gg1$ with some regular assumptions on initial datum and $|\mb{B}_1|>0$. Moreover, the optimal decay rates for $n\geqslant 3$ have been improved by $t^{-\frac{1}{2}}$ if $|\mb{B}_0|=0$, and some different asymptotic profiles appear.
\end{remark}

\section{Large-time behaviors of solutions in the $L^2$ framework}\label{Section-Fourier-space}
We mainly investigate asymptotic behaviors of the solutions $(\rho,v)$ to the linearized compressible Navier-Stokes equations \eqref{Eq_Linearized_NS} in the $L^2$ framework with $\dot{H}^1_1$ datum as $t\gg1$. To begin with, some estimates of solutions in the Fourier space are derived by employing the WKB analysis. Then, we discuss optimal estimates in the sense of same large-time behaviors for the upper and lower bounds under different assumptions on initial datum.
\subsection{Representation of solutions in the Fourier space}
Let us directly apply the partial Fourier transform with respect to spatial variables for the linearized Cauchy problem \eqref{Eq_Linearized_NS} to obtain 
\begin{align}\label{Fourier-Problem}
\begin{cases}
\widehat{\rho}_t+i\gamma\xi\circ\widehat{v}=0,&\xi\in\mb{R}^n, \ t>0,\\
\widehat{v}_t+\alpha|\xi|^2\widehat{v}+\beta\xi(\xi\circ\widehat{v})+i\gamma\xi\widehat{\rho}=0,&\xi\in\mb{R}^n, \ t>0,\\
\widehat{\rho}(0,\xi)=\widehat{\rho}_0(\xi),\ \widehat{v}(0,\xi)=\widehat{v}_0(\xi),&\xi\in\mb{R}^n.
\end{cases}
\end{align}
We next introduce the following zones in the Fourier space:
\begin{align*}
\ml{Z}_{\intt}(\varepsilon_0):=\{|\xi|\leqslant\varepsilon_0\ll1\}, \ \ 
\ml{Z}_{\bdd}(\varepsilon_0,N_0):=\{\varepsilon_0\leqslant |\xi|\leqslant N_0\},\ \  
\ml{Z}_{\extt}(N_0):=\{|\xi|\geqslant N_0\gg1\}.
\end{align*}
Moreover, the cut-off functions $\chi_{\intt}(\xi),\chi_{\bdd}(\xi),\chi_{\extt}(\xi)\in \mathcal{C}^{\infty}$ having supports in their corresponding zones $\ml{Z}_{\intt}(\varepsilon_0)$, $\ml{Z}_{\bdd}(\varepsilon_0/2,2N_0)$ and $\ml{Z}_{\extt}(N_0)$, respectively, such that $\chi_{\bdd}(\xi)=1-\chi_{\intt}(\xi)-\chi_{\extt}(\xi)$ for all $\xi \in \mb{R}^n$. From elementary computations in \cite[Section 2]{Kobayashi-Shibata=2002}, i.e. applying the divergence operator to \eqref{Fourier-Problem}$_2$ associated with \eqref{Fourier-Problem}$_1$, the density solves the viscoelastic damped waves
\begin{align*}
\widehat{\rho}_{tt}+\gamma^2|\xi|^2\widehat{\rho}+(\alpha+\beta)|\xi|^2\widehat{\rho}_t=0.
\end{align*} Then, concerning $\xi\in\ml{Z}_{\intt}(\varepsilon_0)\cup\ml{Z}_{\extt}(N_0)$, the solutions $(\widehat{\rho},\widehat{v})$ to the system \eqref{Fourier-Problem} have the following representations: 
\begin{align}
\widehat{\rho}&=\frac{\lambda_+\mathrm{e}^{\lambda_-t}-\lambda_-\mathrm{e}^{\lambda_+t}}{\lambda_+-\lambda_-}\widehat{\rho}_0-i\gamma\frac{\mathrm{e}^{\lambda_+t}-\mathrm{e}^{\lambda_-t}}{\lambda_+-\lambda_-}\xi\circ\widehat{v}_0,\label{Rep-01}\\
\widehat{v}&=\mathrm{e}^{-\alpha|\xi|^2t}\widehat{v}_0-i\gamma\frac{\mathrm{e}^{\lambda_+t}-\mathrm{e}^{\lambda_-t}}{\lambda_+-\lambda_-}\xi\widehat{\rho}_0+\left(\frac{\lambda_+\mathrm{e}^{\lambda_+t}-\lambda_-\mathrm{e}^{\lambda_-t}}{\lambda_+-\lambda_-}-\mathrm{e}^{-\alpha|\xi|^2t}\right)\frac{\xi(\xi\circ\widehat{v}_0)}{|\xi|^2},\label{Rep-02}
\end{align} 
with the distinct roots of the characteristic equation $\lambda^2+(\alpha+\beta)|\xi|^2\lambda+\gamma^2|\xi|^2=0$, namely,
\begin{align*}
\lambda_{\pm}=-\frac{\alpha+\beta}{2}|\xi|^2\pm\frac{\alpha+\beta}{2}\sqrt{|\xi|^4-\frac{4\gamma^2}{(\alpha+\beta)^2}|\xi|^2}
\end{align*}
due to $|\xi|\leqslant\varepsilon_0\ll 1$ or $|\xi|\geqslant N_0\gg1$.
\subsection{Refined estimates of solutions in the Fourier space}
First of all, considering bounded frequencies $\xi\in\ml{Z}_{\bdd}(\varepsilon_0,N_0)$, thanks to  negative real parts of the characteristic roots, we immediately claim
\begin{align}\label{BDD-Est}
	\chi_{\bdd}(\xi)\left(\,\left|\widehat{\rho}-\mathrm{e}^{-\frac{\alpha+\beta}{2}|\xi|^2t}\widehat{\rho}_0\right|+\left|\widehat{v}-\mathrm{e}^{-\alpha|\xi|^2t}\widehat{v}_0\right|\,\right)\lesssim\chi_{\bdd}(\xi)\mathrm{e}^{-ct}(|\widehat{\rho}_0|+|\widehat{v}_0|)
\end{align}
for large-time.  Indeed, the characteristic roots for $\xi\in\ml{Z}_{\extt}(N_0)$ own the expansions
\begin{align*}
\lambda_+=-\frac{\gamma^2}{\alpha+\beta}+\ml{O}(|\xi|^{-2})\ \ \mbox{and}\ \ 
\lambda_-=-(\alpha+\beta)|\xi|^2+\ml{O}(1).
\end{align*}
As a consequence, according to the representation of solutions, we are able to get
\begin{align}\label{EXTT-Est}
	\chi_{\extt}(\xi)\left(\,\left|\widehat{\rho}-\mathrm{e}^{-\frac{\alpha+\beta}{2}|\xi|^2t}\widehat{\rho}_0\right|+\left|\widehat{v}-\mathrm{e}^{-\alpha|\xi|^2t}\widehat{v}_0\right|\,\right)\lesssim\chi_{\extt}(\xi)\mathrm{e}^{-ct}(|\widehat{\rho}_0|+|\widehat{v}_0|)
\end{align}
for large-time. It means that exponential decay estimates hold for bounded and large frequencies.

We next study behaviors of solutions for small frequencies in a deep way. The characteristic roots $\lambda_{\pm}$ can be expanded by 
\begin{align*}
	\lambda_{\pm}=\pm i\gamma|\xi|-\frac{\alpha+\beta}{2}|\xi|^2\mp \frac{i(\alpha+\beta)^2}{8\gamma}|\xi|^3+\ml{O}(|\xi|^5)
\end{align*}
for $\xi\in\ml{Z}_{\intt}(\varepsilon_0)$. 
According to the representations \eqref{Rep-01} and \eqref{Rep-02}, we may introduce 
\begin{align*}
\widehat{K}_0&:=\frac{1}{|\xi|}\left(\frac{\lambda_+\mathrm{e}^{\lambda_-t}-\lambda_-\mathrm{e}^{\lambda_+t}}{\lambda_+-\lambda_-}-\mathrm{e}^{-\frac{\alpha+\beta}{2}|\xi|^2t}\right)=\frac{\cos(\lambda_{\mathrm{I}}t)}{|\xi|}\mathrm{e}^{\lambda_{\mathrm{R}}t}-\frac{1}{|\xi|}\mathrm{e}^{-\frac{\alpha+\beta}{2}|\xi|^2t}-\frac{\lambda_{\mathrm{R}}\sin(\lambda_{\mathrm{I}}t)}{\lambda_{\mathrm{I}}|\xi|}\mathrm{e}^{\lambda_{\mathrm{R}}t},\\
\widehat{K}_1&:=-i\gamma\frac{\mathrm{e}^{\lambda_+t}-\mathrm{e}^{\lambda_-t}}{\lambda_+-\lambda_-}\frac{\xi}{|\xi|}=-i\gamma\frac{\sin(\lambda_{\mathrm{I}}t)}{\lambda_{\mathrm{I}}}\mathrm{e}^{\lambda_{\mathrm{R}}t}\frac{\xi}{|\xi|},\\
\widehat{K}_2&:=\frac{1}{|\xi|}\left(\frac{\lambda_+\mathrm{e}^{\lambda_+t}-\lambda_-\mathrm{e}^{\lambda_-t}}{\lambda_+-\lambda_-}-\mathrm{e}^{-\alpha|\xi|^2t}\right)=\frac{\cos(\lambda_{\mathrm{I}}t)}{|\xi|}\mathrm{e}^{\lambda_{\mathrm{R}}t}-\frac{1}{|\xi|}\mathrm{e}^{-\alpha|\xi|^2t}+\frac{\lambda_{\mathrm{R}}\sin(\lambda_{\mathrm{I}}t)}{\lambda_{\mathrm{I}}|\xi|}\mathrm{e}^{\lambda_{\mathrm{R}}t},
\end{align*}
because of $\lambda_{\pm}=\lambda_{\mathrm{R}}\pm i\lambda_{\mathrm{I}}$ with
\begin{align*}
	\lambda_{\mathrm{R}}=-\frac{\alpha+\beta}{2}|\xi|^2+\ml{O}(|\xi|^5)\ \ \mbox{and}\ \ \lambda_{\mathrm{I}}=\gamma|\xi|-\frac{(\alpha+\beta)^2}{8\gamma}|\xi|^3+\ml{O}(|\xi|^5).
\end{align*}
Then, we obtain
\begin{align*}
	\chi_{\intt}(\xi)\left(\widehat{\rho}-\mathrm{e}^{-\frac{\alpha+\beta}{2}|\xi|^2t}\widehat{\rho}_0\right)&=\chi_{\intt}(\xi)\left(\widehat{K}_0|\xi|\widehat{\rho}_0+\widehat{K}_1\circ|\xi|\widehat{v}_0\right),\\
	\chi_{\intt}(\xi)\left(\widehat{v}-\mathrm{e}^{-\alpha|\xi|^2t}\widehat{v}_0\right)&=\chi_{\intt}(\xi)\left(\widehat{K}_1|\xi|\widehat{\rho}_0+\widehat{K}_2\frac{\xi(\xi\circ|\xi|\widehat{v}_0)}{|\xi|^2}\right).
\end{align*}

Let us recall the partial Fourier transforms of the multipliers introduced in Theorem \ref{Thm-Optimal-B0}.
A direct subtraction implies
\begin{align}\label{K0}
	\chi_{\intt}(\xi)\left|\widehat{K}_0-\widehat{J}_0\right|
	&\lesssim\chi_{\intt}(\xi)\mathrm{e}^{-c|\xi|^2t}+\chi_{\intt}(\xi)\frac{1}{|\xi|}\left|\cos(\lambda_{\mathrm{I}}t)-\cos(\gamma|\xi|t)\right|\mathrm{e}^{\lambda_{\mathrm{R}}t}\notag\\
	&\quad+\chi_{\intt}(\xi)\frac{1}{|\xi|}|\cos(\gamma|\xi|t)|\left|\mathrm{e}^{\lambda_{\mathrm{R}}t}-\mathrm{e}^{-\frac{\alpha+\beta}{2}|\xi|^2t}\right|\notag\\
	&\lesssim\chi_{\intt}(\xi)\left(1+|\xi|^2t+|\xi|^4t\right)\mathrm{e}^{-c|\xi|^2t}\lesssim \chi_{\intt}(\xi)\mathrm{e}^{-c|\xi|^2t},
\end{align}
where we employed $\cos(\lambda_{\mathrm{I}}t)-\cos(\gamma|\xi|t)=\ml{O}(|\xi|^3)t$
as well as $\mathrm{e}^{\ml{O}(|\xi|^5)t}-1=\ml{O}(|\xi|^5)t\int_0^1\mathrm{e}^{\ml{O}(|\xi|^5)t\tau}\mathrm{d}\tau$.
Analogously,
 we can derive
\begin{align}\label{K1}
	\chi_{\intt}(\xi)\left(\,\left|\widehat{K}_1-\widehat{J}_1\right|+\left|\widehat{K}_2-\widehat{J}_2\right|\,\right)&\lesssim \chi_{\intt}(\xi)\mathrm{e}^{-c|\xi|^2t}.
\end{align}
In conclusion, we claim
\begin{align}
	\chi_{\intt}(\xi)\left|\left(\widehat{\rho}-\mathrm{e}^{-\frac{\alpha+\beta}{2}|\xi|^2t}\widehat{\rho}_0\right)-\left(\widehat{J}_0|\xi|\widehat{\rho}_0+\widehat{J}_1\circ|\xi|\widehat{v}_0\right)\right|&\lesssim\chi_{\intt}(\xi)\mathrm{e}^{-c|\xi|^2t}\left(|\xi|\,|\widehat{\rho}_0|+|\xi|\,|\widehat{v}_0|\right),\label{Err-01}\\
	\chi_{\intt}(\xi)\left|\left(\widehat{v}-\mathrm{e}^{-\alpha|\xi|^2t}\widehat{v}_0\right)-\left(\widehat{J}_1|\xi|\widehat{\rho}_0+\widehat{J}_2\frac{\xi(\xi\circ|\xi|\widehat{v}_0)}{|\xi|^2}\right)\right|
	&\lesssim\chi_{\intt}(\xi)\mathrm{e}^{-c|\xi|^2t}\left(|\xi|\,|\widehat{\rho}_0|+|\xi|\,|\widehat{v}_0|\right).\label{Err-02}
\end{align}
\subsection{Optimal estimates and asymptotic profiles of solutions with $|\mb{B}_0|>0$}
Before deriving the first main result, let us recall the next optimal large-time estimates for the Fourier multipliers from the recent works \cite{Ikehata=2014,Ikehata-Onodera=2017} and \cite[Proposition 4.1]{Chen-Takeda=2022}, respectively.
 \begin{prop}\label{Prop-Fourier-mul-Takeda}
 	Let $c_j>0$ with $j=0,1,2$. The following estimates hold:
	\begin{align*}
	\left\|\chi_{\intt}(\xi)\frac{1}{|\xi|}\sin(c_1|\xi|t)\mathrm{e}^{-c_2|\xi|^2t}\right\|_{L^2}&\simeq\ml{D}_n(t),\\ 
	\left\|\chi_{\intt}(\xi)\frac{1}{|\xi|}\left(\mathrm{e}^{-c_0|\xi|^2t}-\cos(c_1|\xi|t)\mathrm{e}^{-c_2|\xi|^2t}\right)\right\|_{L^2}&\simeq\ml{D}_n(t),
\end{align*}
for any $n\geqslant 1$ and $t\gg1$, where the time-dependent function $\ml{D}_n(t)$ was defined in \eqref{Decay-fun}.
\end{prop}
From the derived estimates \eqref{BDD-Est}-\eqref{Err-02} associated with the triangle inequality, we arrive at 
\begin{align*}
&\left\|\rho(t,\cdot)-\mathrm{e}^{\frac{\alpha+\beta}{2}\Delta t}\rho_0(\cdot)\right\|_{L^2}+\left\|v(t,\cdot)-\mathrm{e}^{\alpha\Delta t}v_0(\cdot)\right\|_{(L^2)^{n}}\\
&\qquad\lesssim\left\|\chi_{\intt}(\xi)\mathrm{e}^{-c|\xi|^2t}\left(|\xi|\,|\widehat{\rho}_0|+|\xi|\,|\widehat{v}_0|\right)\right\|_{L^2}+\left\|\chi_{\intt}(\xi)\left(\widehat{J}_0|\xi|\widehat{\rho}_0+\widehat{J}_1\circ|\xi|\widehat{v}_0\right)\right\|_{L^2}\\
&\qquad\quad+\left\|\chi_{\intt}(\xi)\left(\widehat{J}_1|\xi|\widehat{\rho}_0+\widehat{J}_2\frac{\xi(\xi\circ|\xi|\widehat{v}_0)}{|\xi|^2}\right)\right\|_{(L^2)^n}+\mathrm{e}^{-ct}\left\|\big(\chi_{\bdd}(\xi)+\chi_{\extt}(\xi)\big)\left(|\widehat{\rho}_0|+|\widehat{v}_0|\right)\right\|_{L^2}\\
&\qquad\lesssim\left(\left\|\chi_{\intt}(\xi)\mathrm{e}^{-c|\xi|^2t}\right\|_{L^2}+\sum\limits_{k=0,1,2}\left\|\chi_{\intt}(\xi)\widehat{J}_k\right\|_{L^2}\right)\|(|\xi|\widehat{\rho}_0,|\xi|\widehat{v}_0)\|_{(L^{\infty})^{1+n}}+\mathrm{e}^{-ct}\|(\widehat{\rho}_0,\widehat{v}_0)\|_{(L^2)^{1+n}}\\
&\qquad\lesssim\ml{D}_n(t)\|(\rho_0,v_0)\|_{(\dot{H}^1_1)^{1+n}}+\mathrm{e}^{-ct}\|(\rho_0,v_0)\|_{(L^2)^{1+n}},
\end{align*}
where we used the upper bound estimates in Proposition \ref{Prop-Fourier-mul-Takeda}, as well as the $L^{\infty}$ embedding theorem in the Fourier space. Taking a sufficiently small constant $\alpha_0>0$, we notice that
\begin{align*}
\left|J_0(t,D)|D|\rho_0(x)-J_0(t,x)P_{|D|\rho_0}\right|&\leqslant\left(\int_{|y|\leqslant t^{\alpha_0}}+\int_{|y|\geqslant t^{\alpha_0}}\right)|J_0(t,x-y)-J_0(t,x)| \big||D|\rho_0(y)\big|\mathrm{d}y\\
&\lesssim t^{\alpha_0}\int_{|y|\leqslant t^{\alpha_0}}|\nabla J_0(t,x-\theta_0y)|\,\big||D|\rho_0(y)\big|\mathrm{d}y\\
&\quad+\int_{|y|\geqslant t^{\alpha_0}}\big(|J_0(t,x-y)|+|J_0(t,x)|\big)\big| |D|\rho_0(y)\big|\mathrm{d}y,
\end{align*}
in which the next estimate was applied:
\begin{align*}
	|J_0(t,x-y)-J_0(t,x)|\lesssim|y|\,|\nabla J_0(t,x-\theta_0y)|
\end{align*}
with $\theta_0\in(0,1)$. Hence, it leads to
\begin{align}
&\left\|\chi_{\intt}(D)\left(J_0(t,D)|D|\rho_0(\cdot)-J_0(t,\cdot)P_{|D|\rho_0}\right)\right\|_{L^2}\notag\\
&\qquad\lesssim t^{\alpha_0}\left\|\chi_{\intt}(\xi)|\xi|\widehat{J}_0\right\|_{L^2}\|\rho_0\|_{\dot{H}^1_1}+\left\|\chi_{\intt}(\xi)\widehat{J}_0\right\|_{L^2}\int_{|y|\geqslant t^{\alpha_0}}\big||D|\rho_0(y)\big|\mathrm{d}y=o\big(\ml{D}_n(t)\big)\label{Est_001}
\end{align}
as $t\gg1$ by choosing small $\alpha_0>0$ and $|D|\rho_0\in L^1$. With the same approach, we obtain
\begin{align*}
\left\|\chi_{\intt}(D)\left[\big(J_0(t,D)|D|\rho_0(\cdot)+J_1(t,D)\circ|D|v_0(\cdot)\big)-\big(J_0(t,\cdot)P_{|D|\rho_0}+J_1(t,\cdot)\circ P_{|D|v_0}\big)\right]\right\|_{L^2}=o\big(\ml{D}_n(t)\big)
\end{align*}
as large-time $t\gg1$. Let us recall the facts that
\begin{align*}
\int_{\mb{S}^{n-1}}\omega_j\omega_k\mathrm{d}\sigma_{\omega}=\begin{cases}
\displaystyle{\frac{1}{n}|\mb{S}^{n-1}|}&\mbox{when}\ \ j=k,\\
0&\mbox{when}\ \ j\neq k,
\end{cases}
\end{align*}
and the expression with a scalar function $\widehat{g}(t,|\xi|)$ as follows:
\begin{align}
\left\|\widehat{g}(t,|\xi|)\frac{\xi}{|\xi|}\circ P_{|D|v_0}\right\|_{L^2}^2&=\int_{\mb{R}^n}|\widehat{g}(t,|\xi|)|^2\frac{1}{|\xi|^2}\left(\xi_1P_{|D|v_0^{(1)}}+\cdots+\xi_nP_{|D|v_0^{(n)}}\right)^2\mathrm{d}\xi\notag\\
&=\frac{|\mb{S}^{n-1}|}{n}\int_0^{\infty}|\widehat{g}(t,r)|^2r^{n-1}\mathrm{d}r\,|P_{|D|v_0}|^2.\label{non-radial}
\end{align} Finally, concerning large-time $t\gg1$, with the aid of Minkowski's inequality, one may see
\begin{align*}
&\left\|\rho(t,\cdot)-\mathrm{e}^{\frac{\alpha+\beta}{2}\Delta t}\rho_0(\cdot)\right\|_{L^2}\geqslant\left\|\chi_{\intt}(D)\left(\rho(t,\cdot)-\mathrm{e}^{\frac{\alpha+\beta}{2}\Delta t}\rho_0(\cdot)\right)\right\|_{L^2}\\
&\qquad\geqslant\left\|\chi_{\intt}(D)\big(J_0(t,\cdot)P_{|D|\rho_0}+J_1(t,\cdot)\circ P_{|D|v_0}\big)\right\|_{L^2}\\
&\qquad\quad-\left\|\chi_{\intt}(D)\left[\big(J_0(t,D)|D|\rho_0(\cdot)+J_1(t,D)\circ|D|v_0(\cdot)\big)-\big(J_0(t,\cdot)P_{|D|\rho_0}+J_1(t,\cdot)\circ P_{|D|v_0}\big)\right]\right\|_{L^2}\\
&\qquad\quad-\left\|\chi_{\intt}(D)\left[\left(\rho(t,\cdot)-\mathrm{e}^{\frac{\alpha+\beta}{2}\Delta t}\rho_0(\cdot)\right)-\big(J_0(t,D)|D|\rho_0(\cdot)+J_1(t,D)\circ|D|v_0(\cdot)\big) \right]\right\|_{L^2}\\
&\qquad\gtrsim\left\|\chi_{\intt}(D)\big(J_0(t,\cdot)P_{|D|\rho_0}+J_1(t,\cdot)\circ P_{|D|v_0}\big)\right\|_{L^2}-o\big(\ml{D}_n(t)\big)-t^{-\frac{n}{4}}\|(\rho_0,v_0)\|_{(\dot{H}^1)^{1+n}}.
\end{align*}
Moreover, by applying \eqref{non-radial} and Proposition \ref{Prop-Fourier-mul-Takeda} from the lower bounds, we have
\begin{align*}
&\left\|\chi_{\intt}(D)\big(J_0(t,\cdot)P_{|D|\rho_0}+J_1(t,\cdot)\circ P_{|D|v_0}\big)\right\|_{L^2}^2\\
&\qquad\gtrsim\left\|\chi_{\intt}(\xi)\frac{1}{|\xi|}\big(\cos(\gamma|\xi|t)-1\big)\mathrm{e}^{-\frac{\alpha+\beta}{2}|\xi|^2t}\right\|_{L^2}^2|P_{|D|\rho_0}|^2+\left\|\chi_{\intt}(\xi)\frac{\sin(\gamma|\xi|t)}{|\xi|}\mathrm{e}^{-\frac{\alpha+\beta}{2}|\xi|^2t}\frac{\xi}{|\xi|}\circ P_{|D|v_0}\right\|_{L^2}^2\\
&\qquad\gtrsim \big(\ml{D}_n(t)\big)^2\left(|P_{|D|\rho_0}|^2+|P_{|D|v_0}|^2\right),
\end{align*}
which shows immediately (see the definition \eqref{Intro_B0} of $|\mb{B}_0|^2$)
\begin{align*}
	\left\|\rho(t,\cdot)-\mathrm{e}^{\frac{\alpha+\beta}{2}\Delta t}\rho_0(\cdot)\right\|_{L^2}\gtrsim \ml{D}_n(t)|\mb{B}_0|-o\big(\ml{D}_n(t)\big)-t^{-\frac{n}{4}}\|(\rho_0,v_0)\|_{(\dot{H}^1_1)^{1+n}}\gtrsim \ml{D}_n(t)|\mb{B}_0|.
\end{align*}
Similarly to the above computations, the lower bound for another quantity also can be estimated
\begin{align*}
	\left\|v(t,\cdot)-\mathrm{e}^{\alpha\Delta t}v_0(\cdot)\right\|_{(L^2)^{n}}\gtrsim \ml{D}_n(t)|\mb{B}_0|
\end{align*}
for large-time $t\gg1$. Particularly, focusing on the lower-dimensional cases $n=1,2$, due to the boundedness
\begin{align*}
\left\|\mathrm{e}^{\frac{\alpha+\beta}{2}\Delta t}\rho_0(\cdot)\right\|_{L^2}+\left\|\mathrm{e}^{\alpha\Delta t}v_0(\cdot)\right\|_{(L^2)^n}\leqslant\|(\rho_0,v_0)\|_{(L^2)^{1+n}},
\end{align*} 
the optimal estimates can be derived immediately
\begin{align*}
\left\|(\rho,v)(t,\cdot)\right\|_{(L^2)^{1+n}}&\lesssim \ml{D}_n(t)\|(\rho_0,v_0)\|_{(\dot{H}^1_1)^{1+n}}+\|(\rho_0,v_0)\|_{(L^2)^{1+n}}\lesssim \ml{D}_n(t)\|(\rho_0,v_0)\|_{(L^2\cap \dot{H}^{1}_1)^{1+n}},\\
\left\|(\rho,v)(t,\cdot)\right\|_{(L^2)^{1+n}}&\gtrsim\ml{D}_n(t)|\mb{B}_0|-\|(\rho_0,v_0)\|_{(L^2)^{1+n}}\gtrsim\ml{D}_n(t)|\mb{B}_0|,
\end{align*}
as $t\gg1$ and $n=1,2$. The desired estimates of solutions in Theorem \ref{Thm-Optimal-B0} are finished.

We now turn to asymptotic profiles of the solutions. Let us recall \cite[Lemma 2.2]{Ikehata=2014}, namely,
\begin{align}\label{G}
	|\widehat{g}-P_g|\lesssim|\xi|\,\|g\|_{L^{1,1}}.
\end{align}
By assuming additionally $(|D|\rho_0,|D|v_0)\in (L^{1,1})^{1+n}$, the error estimates
\begin{align*}
&\left\|\chi_{\intt}(D)\left[\big(J_0(t,D)|D|\rho_0(\cdot)+J_1(t,D)\circ|D|v_0(\cdot)\big)-\big(J_0(t,\cdot)P_{|D|\rho_0}+J_1(t,\cdot)\circ P_{|D|v_0}\big)\right]\right\|_{L^2}\\
&\qquad\lesssim\|\chi_{\intt}(\xi)|\xi|\widehat{J}_0\|_{L^2}\|\,|D|\rho_0\|_{L^{1,1}}+\|\chi_{\intt}(\xi)|\xi|\widehat{J}_1\|_{(L^2)^n}\|\,|D|v_0\|_{(L^{1,1})^n}\\
&\qquad\lesssim t^{-\frac{n}{4}}\|(|D|\rho_0,|D|v_0)\|_{(L^{1,1})^{1+n}}
\end{align*}
for large-time $t\gg1$, improve the previous one \eqref{Est_001}.
 Applying these error estimates as well as
\begin{align*}
	\left\|\big(1-\chi_{\intt}(D)\big)\big(J_0(t,\cdot)P_{|D|\rho_0}+J_1(t,\cdot)\circ P_{|D|v_0}\big)\right\|_{L^2}\lesssim \mathrm{e}^{-ct}\left(|P_{|D|\rho_0}|+|P_{|D|v_0}|\right),
\end{align*}
from the triangle inequality, we may derive
\begin{align*}
&\left\|\rho(t,\cdot)-\mathrm{e}^{\frac{\alpha+\beta}{2}\Delta t}\rho_0(\cdot)-\big(J_0(t,\cdot)P_{|D|\rho_0}+J_1(t,\cdot)\circ P_{|D|v_0}\big)\right\|_{L^2}\\
&\qquad\lesssim t^{-\frac{n}{4}}\|(|D|\rho_0,|D|v_0)\|_{(L^{1,1})^{1+n}}+\mathrm{e}^{-ct}\|(\rho_0,v_0)\|_{(L^2)^{1+n}}
\end{align*}
for $t\gg1$. By the same way, the refined estimates \eqref{Est-002} also can be obtained. The proof of Theorem \ref{Thm-Optimal-B0} is completed.

\subsection{Optimal estimates and asymptotic profiles of solutions with $|\mb{B}_0|=0$}
Indeed, our current consideration $|\mb{B}_0|=0$ is equivalent to $|P_{|D|\rho_0}|=0$ as well as $|P_{|D|v_0}|=0$. Therefore, by the same approach of proving \eqref{Est-002}, the upper bounds estimates can be obtained
\begin{align*}
\left\|\left(\rho-\mathrm{e}^{\frac{\alpha+\beta}{2}\Delta t}\rho_0\right)(t,\cdot)\right\|_{L^2}+\left\|\left(v-\mathrm{e}^{\alpha\Delta t}v_0\right)(t,\cdot)\right\|_{(L^2)^n}\lesssim t^{-\frac{n}{4}}\left(\|(|D|\rho_0,|D|v_0)\|_{(L^{1,1})^{1+n}}+\|(\rho_0,v_0)\|_{(L^2)^{1+n}}\right),
\end{align*}
for $t\gg1$. Moreover, we may derive
\begin{align*}
\left\|\mathrm{e}^{\frac{\alpha+\beta}{2}\Delta t}\rho_0(\cdot)\right\|_{L^2}+\left\|\mathrm{e}^{\alpha\Delta t}v_0(\cdot)\right\|_{(L^2)^n}&\lesssim\left\|\chi_{\intt}(\xi)\frac{\mathrm{e}^{-c|\xi|^2t}}{|\xi|}\left(|\xi|\,\|\,|D|\rho_0\|_{L^{1,1}}+|\xi|\,\|\,|D|v_0\|_{(L^{1,1})^n}\right)\right\|_{L^2} \\
&\quad+\left\|\big(1-\chi_{\intt}(\xi)\big)\mathrm{e}^{-c|\xi|^2t}\left(|\widehat{\rho}_0|+|\widehat{v}_0|\right)\right\|_{L^2}\\
&\lesssim t^{-\frac{n}{4}}\left(\|(|D|\rho_0,|D|v_0)\|_{(L^{1,1})^{1+n}}+\|(\rho_0,v_0)\|_{(L^2)^{1+n}}\right),
\end{align*}
where we used \eqref{G} with $P_g=0$ in this case. By employing the triangle inequality,
it results the upper bounds estimates of the solutions in Theorem \ref{Thm_Optimal-B1}.

By taking a further assumption on initial datum such that $(|D|\rho_0,|D|v_0)\in (L^{1,2})^{1+n}$, the next crucial expansion:
\begin{align}\label{Expan}
\widehat{g}=P_g+i\xi\circ Q_g+\widehat{E}_{g}
\end{align}
holds with the estimate (see \cite[Lemma 5.1]{Ikehata-Michihisa=2019} in detail)
\begin{align}\label{GG}
|\widehat{E}_g|\lesssim|\xi|^2\|g\|_{L^{1,2}}.
\end{align}
Note that the expansion \eqref{Expan} is helpful for us to understand \eqref{G}.
Therefore, plugging $g=|D|\rho_0$ in \eqref{Expan}, the solution $\widehat{\rho}$ for small frequencies can be represented by
\begin{align*}
\chi_{\intt}(\xi)\widehat{\rho}&=\chi_{\intt}(\xi)\left(\frac{1}{|\xi|}\mathrm{e}^{-\frac{\alpha+\beta}{2}|\xi|^2t}+\widehat{K}_0\right)\left(i\xi\circ Q_{|D|\rho_0}+\widehat{E}_{|D|\rho_0}\right)+\chi_{\intt}(\xi)\sum\limits_{k=1}^{n}\widehat{K}_1^{(k)}\left(i\xi\circ Q_{|D|v_0^{(k)}}+\widehat{E}_{|D|v_0^{(k)}}\right)
\end{align*}
with $\widehat{K}_1=(\widehat{K}_1^{(1)},\dots,\widehat{K}_1^{(n)})$. For this reason,  we know
\begin{align*}
&\chi_{\intt}(\xi)\left[\widehat{\rho}-\left(\frac{1}{|\xi|}\mathrm{e}^{-\frac{\alpha+\beta}{2}|\xi|^2t}+\widehat{J}_0\right)i\xi\circ Q_{|D|\rho_0}-\sum\limits_{k=1}^n\widehat{J}_1^{(k)}i\xi\circ Q_{|D|v_0^{(k)}}\right]\\
&\qquad=\chi_{\intt}(\xi)\left[\frac{1}{|\xi|}\widehat{E}_{|D|\rho_0}\mathrm{e}^{-\frac{\alpha+\beta}{2}|\xi|^2t}+\left(\widehat{K}_0-\widehat{J}_0\right)i\xi\circ Q_{|D|\rho_0}+\widehat{K}_0\widehat{E}_{|D|\rho_0}\right]\\
&\qquad\quad+\chi_{\intt}(\xi)\sum\limits_{k=1}^n\left[\left(\widehat{K}_1^{(k)}-\widehat{J}_1^{(k)}\right)i\xi\circ Q_{|D|v_0^{(k)}}+\widehat{K}_1^{(k)}\widehat{E}_{|D|v_0^{(k)}}\right].
\end{align*}
Recalling the definition of $\widetilde{\rho}^0(t,x)$ in Theorem \ref{Thm_Optimal-B1}, combining with \eqref{K0}, \eqref{K1} and \eqref{GG}, the next error estimates:
\begin{align}\label{two-star}
\left\|\left(\rho-\widetilde{\rho}^0\right)(t,\cdot)\right\|_{L^2}&\lesssim\left\|\chi_{\intt}(\xi)|\xi|\mathrm{e}^{-c|\xi|^2t}\right\|_{L^2}\|(|D|\rho_0,|D|v_0)\|_{(L^{1,2})^{1+n}}+\mathrm{e}^{-ct}\|(\rho_0,v_0)\|_{(L^2)^{1+n}}\notag\\
&\lesssim t^{-\frac{1}{2}-\frac{n}{4}}\left(\|(|D|\rho_0,|D|v_0)\|_{(L^{1,2})^{1+n}}+\|(\rho_0,v_0)\|_{(L^2)^{1+n}}\right)
\end{align}
hold for $t\gg1$, where we applied
\begin{align*}
&\left\|\big(1-\chi_{\extt}(\xi)\big)\left[\left(\frac{1}{|\xi|}\mathrm{e}^{-\frac{\alpha+\beta}{2}|\xi|^2t}+\widehat{J}_0\right)i\xi\circ Q_{|D|\rho_0}+\sum\limits_{k=1}^n\widehat{J}^{(k)}_1i\xi\circ Q_{|D|v_0^{(k)}} \right]\right\|_{L^2}\\
&\qquad\lesssim \mathrm{e}^{-ct}\left(|Q_{|D|\rho_0}|+\sum\limits_{k=1}^n|Q_{|D|v_0^{(k)}}|\right).
\end{align*}
 Similarly, due to 
\begin{align*}
&\chi_{\intt}(\xi)\left|\widehat{v}-\widehat{J}_1i\xi\circ Q_{|D|\rho_0}-\frac{i\xi}{|\xi|}\mathrm{e}^{-\alpha|\xi|^2t}\circ Q_{|D|v_0}-\widehat{J}_2\frac{\xi}{|\xi|^2}\sum\limits_{k=1}^n\xi_ki\xi\circ Q_{|D|v_0^{(k)}}\right|\\
&\qquad\lesssim \chi_{\intt}(\xi)|\xi|\mathrm{e}^{-c|\xi|^2}\|(|D|\rho_0,|D|v_0)\|_{(L^{1,2})^{1+n}},
\end{align*}
we can proved the error estimates
\begin{align*}
\left\|\left(v-\widetilde{v}^0\right)(t,\cdot)\right\|_{(L^2)^n}\lesssim t^{-\frac{1}{2}-\frac{n}{4}}\left(\|(|D|\rho_0,|D|v_0)\|_{(L^{1,2})^{1+n}}+\|(\rho_0,v_0)\|_{(L^2)^{1+n}}\right)
\end{align*} for large-time $t\gg1$.

Finally, we turn to the lower bounds estimates. According to the definition of $\widehat{J}_0$ and $\widehat{J}_1$, we may express the profile in the $L^2$ norm as follows:
\begin{align*}
\left\|\widetilde{\rho}^0(t,\cdot)\right\|_{L^2}^2&=\left\|\frac{\cos(\gamma|\xi|t)}{|\xi|}\mathrm{e}^{-\frac{\alpha+\beta}{2}|\xi|^2t}\xi\circ Q_{|D|\rho_0}\right\|_{L^2}^2+\left\|\sum\limits_{k=1}^n\frac{\sin(\gamma|\xi|t)}{|\xi|}\mathrm{e}^{-\frac{\alpha+\beta}{2}|\xi|^2t}\xi_k\frac{\xi}{|\xi|}\circ Q_{|D|v_0^{(k)}}\right\|_{L^2}^2\\
&=\frac{|\mb{S}^{n-1}|}{n}\int_0^{\infty}|\cos(\gamma rt)|^2r^{n-1}\mathrm{e}^{-(\alpha+\beta)r^2t}\mathrm{d}r\,|Q_{|D|\rho_0}|^2\\
&\quad+\int_0^{\infty}|\sin(\gamma rt)|^2r^{n-1}\mathrm{e}^{-(\alpha+\beta)r^2t}\mathrm{d}r\sum\limits_{k=1}^n\int_{\mb{S}^{n-1}}\omega_k^2\left(\omega\circ Q_{|D|v_0^{(k)}}\right)^2\mathrm{d}\sigma_{\omega}.
\end{align*}
Clearly, by using the Riemann-Lebesgue theorem associated with $2\sin^2z=1-\cos(2z)$ and $2\cos^2z=1+\cos(2z)$, one may follow the similar methods to the one in \cite{Ikehata=2014}, and get
\begin{align*}
c_0\int_0^{\infty}|\cos(\gamma rt)|^2r^{n-1}\mathrm{e}^{-(\alpha+\beta)r^2t}\mathrm{d}r+c_1\int_0^{\infty}|\sin(\gamma rt)|^2r^{n-1}\mathrm{e}^{-(\alpha+\beta)r^2t}\mathrm{d}r\gtrsim t^{-\frac{n}{2}} 
\end{align*}
as $t\gg1$ with $c_0,c_1>0$. Let us recall the relation (see, for example, \cite[Corollary 2.7]{Takeda=2022})
\begin{align*}
\int_{\mb{S}^{n-1}}\omega_k^2(\omega\circ Q_g)^2\mathrm{d}\sigma_{\omega}=\frac{\pi^{\frac{n}{2}}}{\Gamma(\frac{n+4}{2})}\left(\frac{|Q_g|^2}{2}+|Q_g^{(k)}|^2\right),
\end{align*}
where $\Gamma(s)=\int_0^{\infty}\mathrm{e}^{-z}z^{s-1}\mathrm{d}z$ is the Gamma function. Thus, the lower bound of the profile is
\begin{align*}
\left\|\widetilde{\rho}^0(t,\cdot)\right\|_{L^2}^2\gtrsim t^{-\frac{n}{2}}\left(|Q_{|D|\rho_0}|^2+\sum\limits_{k=1}^n\left(|Q_{|D|v_0^{(k)}}|^2+|Q_{|D|v_0^{(k)}}^{(k)}|^2\right) \right).
\end{align*}
From \eqref{two-star}, one derives
\begin{align*}
\|\rho(t,\cdot)\|_{L^2}\gtrsim t^{-\frac{n}{4}}|\mb{B}_1|-t^{-\frac{1}{2}-\frac{n}{4}}\left(\|(|D|\rho_0,|D|v_0)\|_{(L^{1,2})^{1+n}}+\|(\rho_0,v_0)\|_{(L^2)^{1+n}}\right)\gtrsim t^{-\frac{n}{4}}|\mb{B}_1|
\end{align*}
as large-time. Then, we can derive the lower bounds estimates for $\|v(t,\cdot)\|_{(L^2)^n}$ in an analogous way. All in all, the proof of Theorem \ref{Thm_Optimal-B1} is completed.

\section{Final remarks}\label{Section_Final_Remark}
Throughout  this paper, we have succeeded in deriving not only new upper bounds estimates but also optimal lower bounds estimates of the solutions $(\rho,v)$ to the linearized compressible Navier-Stokes equations \eqref{Eq_Linearized_NS} with initial datum $(|D|\rho_0,|D|v_0)$ belonging to some weighted $L^1$ spaces. Moreover, we introduce a new threshold to distinguish two different large-time behaviors of the solutions. Clearly, this is to state that the obtained large-time estimates in this work are actually sharp. More generally, we expect that our philosophy utilized in this paper can be applied to study optimal estimates associated with some new thresholds for the linearized Cauchy problem in compressible fluids, for example, the linearized Navier-Stokes-Fourier equations \cite{Kobayashi-Shibata=1999}, the linearized Navier-Stokes equations of Korteweg-type \cite{Kawashima-Shibata-Xu=2022}, the linearized Navier-Stokes equations with capillarity \cite{Danchin-Des=2001,Kawa-Shibata-Xu=2021}, and the linearized Navier-Stokes-Poisson equations \cite{Li-Matsumura-Zhang=2010}.

\section*{Acknowledgments}
The second author was supported in part by Grant-in-Aid for scientific Research (C) 20K03682 of JSPS.

\end{document}